\documentclass[12pt,reqno]{amsart}

\usepackage{hyperref}  %provides clickable links to references within the document

\usepackage{graphicx}
\usepackage{amsfonts}
\usepackage{amssymb}
\usepackage{amsmath}
\usepackage{amsthm}

\newtheorem{theorem}{Theorem}[section]

\newcommand{\bF}{\mathbb F}

\begin{document}

\sloppy

\title{OpenAI's proof of the Cycle Double Cover Theorem}

\author{Jim Geelen}

\date{\today}

\maketitle

\section{Introduction}

The following notes are my attempt to clarify, to myself, the OpenAI proof of the Cycle Double Cover Conjecture~\cite{OpenAI}.
This write-up offers no new insights, but some people may find it easier to read.
\begin{theorem}[Cycle Double Cover Theorem]
Every bridgeless graph has a cycle double cover.
\end{theorem}

An {\em edge-bicolouring} of a graph $G$ is an assignment of pairs of colours to the edges of a graph such that each colour occurs an even number of times at each vertex. Since each colour-class decomposes into edge-disjoint cycles,
the existence of an edge-bicolouring is equivalent to the existence of a cycle double cover.

Jaeger~\cite{Jaeger} proved that every bridgeless graph admits a nowhere-zero $\bF_2^3$-flow.
(A nowhere-zero $\bF_2^3$-flow is constructed by assigning each edge a non-zero ``flow-value" from $\bF_2^3$, such that at each vertex the flow-values of its incident edges sum to zero.)
OpenAI uses linear algebra to construct an edge-bicolouring from a nowhere-zero $\bF_2^3$-flow. 

\section{The linear algebra formulation}

It suffices to prove the Cycle Double Cover Theorem for cubic graphs (this is well-known and easy to prove by ``uncontracting" edges).
Let $f:E\rightarrow \bF_2^3$ be a nowhere-zero $\bF_2^3$-flow on a bridgeless cubic graph $G=(V,E)$. 

Let $a$, $b$ and $c$ be the edges incident with a vertex $v$. We define
$\alpha_{v,a}=f_b$ and $\beta_{v,a}=f_c$.
This definition depends on the ordering of the other two incident edges, but we can define these values across the graph in such a way that
$$\{\alpha_{v,a},\alpha_{v,b},\alpha_{v,c}\} = \{f_a,f_b,f_c\}.$$
Note that each element of $\bF_2^3$ occurs an even number of times among
$$ (\alpha_{v,a},\alpha_{v,b},\alpha_{v,c},\beta_{v,a},\beta_{v,b},\beta_{v,c}).$$

Our goal is to find ``vertex potentials" $p:V\rightarrow \bF_2^3$ so that for each edge $e=uv$ we have
\begin{equation}\label{equal}
 \{p_v+\alpha_{v,e},p_v+\beta_{v,e}\} = \{p_u+\alpha_{u,e},p_u+\beta_{u,e}\}.
\end{equation}
The existence of such potentials will give us an edge-bicolouring and hence a cycle double cover.

Since the flow-values on any cut sum to zero, for each edge $e=uv$,
$$ \alpha_{v,e}+\beta_{v,e} = \alpha_{u,e}+\beta_{u,e}.$$
But then 
$$ p_v+\alpha_{v,e}=p_u+\alpha_{u,e} \quad \Leftrightarrow\quad p_v+\beta_{v,e}=p_u+\beta_{u,e},\mbox{ and}$$
$$ p_v+\alpha_{v,e}=p_u+\beta_{u,e} \quad \Leftrightarrow\quad p_v+\beta_{v,e}=p_u+\alpha_{u,e}.$$
Thus, condition~(\ref{equal}) simplifies to
$$ p_v+\alpha_{v,e} \in \{p_u+\alpha_{u,e},p_u+\beta_{u,e}\}.$$
By defining a variable $z_e\in\bF_2$, we can rewrite this as
$$ p_v+\alpha_{v,e} = p_u+\alpha_{u,e} + z_e(\alpha_{u,e}+\beta_{u,e}).$$
Moreover, we have $\alpha_{u,e}+\beta_{u,e}=f_e.$

In summary, there exists a feasible potential if and only if there exist $p:V\rightarrow \bF_2^3$ and $z:E\rightarrow \bF_2$ such that
for each edge $e=uv$ we have
\begin{equation}\label{formulation}
 p_u+p_v +f_ez_e = \alpha_{u,e}+\alpha_{v,e}.
\end{equation}
Each equation in~(\ref{formulation}) can be separated into three equations over $\bF_2$.
Therefore~(\ref{formulation}) is a system of linear equations over $\bF_2$.

\section{Completing the proof}

We proceed by way of contradiction. Suppose that $G$ has no cycle double cover and hence~(\ref{formulation}) has no solution.
If a system of linear equations over $\bF_2$ has no solution, then some subset of those equations adds to the absurd equation $0=1$.
So there exists $y:E\rightarrow \bF_2^3$ such that
\begin{equation}\label{constantterm}
\sum_{e=uv\in E} (\alpha_{u,e}+\alpha_{v,e})\cdot y_e = 1, \mbox{ and}
\end{equation}
$$ \sum_{e=uv\in E} (p_u+p_v +f_ez_e)\cdot y_e = 0.$$
We rewrite the second equation as
\begin{equation}\label{variableterm}
 \sum_{v\in V} p_v\cdot\left(\sum_{e\sim v} y_e\right) + \sum_{e\in E} (f_e\cdot y_e) z_e =0. 
\end{equation}
Here $e\sim v$ denotes incidence. The coefficient of each variable in~(\ref{variableterm}) is zero, so
for each vertex $v$,
\begin{equation}\label{flow}
 \sum_{e\sim v} y_e  =0,
\end{equation}
and, for each edge $e$,
\begin{equation}\label{orthogonal}
  f_e\cdot y_e=0.
\end{equation}

We claim that for each vertex $v$ there exists $\delta_v\in\bF_2$ such that for each edge $e\sim v$ we have
$\alpha_{v,e}\cdot y_e = \delta_v$. In order to prove this let $a$, $b$ and $c$ be the edges incident with $v$.
Up to symmetry we may assume that $\alpha_{v,a}=f_b$ and $\alpha_{v,b}=f_c$ and $\alpha_{v,c}=f_a$. Now, by~(\ref{flow})
and~(\ref{orthogonal}),
\begin{eqnarray*}
\alpha_{v,a}\cdot y_a + \alpha_{v,b} \cdot y_b &=&  f_b\cdot y_a + f_c \cdot y_b\\
&=& (f_a+f_c)\cdot y_a + f_c \cdot (y_a+y_c)\\
&=& 0.
\end{eqnarray*}
Thus $\alpha_{v,a}\cdot y_a = \alpha_{v,b} \cdot y_b$ and that proves the claim. 

We further claim that, if $\delta_v=0$,
then either one or all three of the vectors $y_a$, $y_b$, and $y_c$ is zero. When $\delta_v=0$, we have 
$f_a\cdot y_a =0$ and $f_b \cdot y_a = \alpha_{v,a}\cdot y_a = \delta_v=0$ and
$f_c \cdot y_a = (f_a+f_b)\cdot y_a =0$. Similarly both $y_b$ and $y_c$ are orthogonal to each of
$f_a$, $f_b$, and $f_c$. There is only one non-zero vector in $\bF^3_2$ that is orthogonal to each of 
$f_a$, $f_b$, and $f_c$. Then, since $y_a+y_b+y_c=0$, one of $y_a$, $y_b$, and $y_c$ is zero and the other
two vectors are the same (possibly also zero), proving the second claim.

To complete the proof, let $G^+$ denote the subgraph of $G$ obtained by deleting any edge $e$ with $y_e=0$ and let 
$X$ be the set of vertices $v$ with $\delta_v=0$. Note that each vertex in $X$ has even degree in $G^+$,
so there are an even number of edges in $G^+$ with exactly one end in $X$.

Consider an edge $e=uv$. Note that $(\alpha_{v,e}+\alpha_{u,e})\cdot y_e = \delta_u +\delta_v$ which is $1$
if and only if exactly one of $u$ and $v$ is in $X$. Moreover, if $u\in X$ and $v\not\in X$, then
$\alpha_{v,e}\cdot y_e=1$ so $y_e\neq 0$ and hence $e\in E(G^+)$.
But this contradicts~(\ref{constantterm}).

Checkmate!

\end{document}